\documentclass[11pt,a4paper]{amsart}
\usepackage{amsmath,amssymb,amsthm,amsxtra,latexsym}
\newtheorem{theorem}{Theorem}
\newtheorem{lemma}[theorem]{Lemma}
\newtheorem{proposition}[theorem]{Proposition}
\newtheorem{corollary}[theorem]{Corollary}
\theoremstyle{definition}

\theoremstyle{remark}
\newtheorem{remark}{Remark}
\DeclareMathOperator{\sech}{sech}
\newcommand{\R}{\mathbb{R}}
\newcommand{\C}{\mathbb{C}}
\newcommand{\N}{\mathbb{N}}
\newcommand{\Z}{\mathbb{Z}}
\newcommand{\la}{\langle}
\newcommand{\ra}{\rangle}
\newcommand{\pd}{\partial}
\newcommand{\Lr}{L^2_{r}(\R^2)}
\begin{document}
\title{Instability of vortex solitons for 2D focusing NLS}
\thanks{This research is supported by Grant-in-Aid for Scientific
Research (No. 17740079)).}
\author{Tetsu Mizumachi}
\address{Faculty of Mathematics, Kyushu University, 
Hakozaki 6-10-1, 812-8581 Japan.}
\email{mizumati@math.kyushu-u.ac.jp}
\begin{abstract}
We study instability of a vortex soliton
$e^{i(m\theta+\omega t)}\phi_{\omega,m}(r)$ to
$$iu_t+\Delta u+|u|^{p-1}u=0,\quad\text{for $x\in\R^n$, $t>0$,}$$
where $n=2$, $m\in\N$ and  $(r,\theta)$ are polar coordinates in $\R^2$.
Grillakis \cite{Gr} proved that every radially standing wave solutions
are unstable if $p>1+4/n$. However, we do not have any examples of
unstable standing wave solutions in the subcritical case $(p<1+n/4)$.

Suppose $\phi_{\omega,m}$ is nonnegative. 
We investigate a limiting profile of $\phi_{\omega,m}$ as $m\to\infty$
and  prove that for every $p>1$, there exists an $m_*\in \N$ such that
for $m\ge m_*$, a vortex soliton
$e^{i(m\theta+\omega t)}\phi_{\omega,m}(r)$ becomes unstable to the
perturbations of the form $e^{i(m+j)\theta}v(r)$
with $1\ll j\ll m$.
\end{abstract}
\maketitle
\section{Introduction}
\label{intro}
In the present paper, we consider instability of radially symmetric
vortex solitons to 2-dimensional nonlinear Schr\"odinger equations
\begin{equation}
  \label{eq:NLS}
\left\{
  \begin{aligned}
& iu_t+\Delta u+f(u)=0 \quad\text{for 
$(x,t)\in\mathbb{R}^n\times\mathbb{R}$},
\\
& u(x,0)=u_0(x)\quad\text{for $x\in\mathbb{R}^2$,}
  \end{aligned}\right.
\end{equation}
where $n=2$ and $f(u)=|u|^{p-1}u.$
Let $\omega>0$, $m\in\N\cup\{0\}$, and let
$e^{i(\omega t+m\theta)}\phi_\omega(r)$ be a standing
wave solution of \eqref{eq:NLS} belonging to $H^1(\R^2)$.
Here $r$ and $\theta$ denote polar coordinates in $\R^2$.
Then $\phi_\omega(r)$ is a solution to
\begin{equation}
\label{eq:B}
\left\{
  \begin{aligned}
& \phi''+\frac{1}{r}\phi'-\left(\omega+\frac{m^2}{r^2}\right)\phi+f(\phi)=0
\quad\text{for $r>0$},
\\
& \lim_{r\to0}\frac{\phi(r)}{r^m}=\lim_{r\to0}\frac{\phi'(r)}{mr^{m-1}},\\
& \lim_{r\to\infty}\phi(r)=0.
  \end{aligned}\right.
\end{equation}
We remark that
$e^{im\theta}\phi_\omega(r)$ is a solution to the scalar field equation
\begin{equation}
  \label{eq:sfe}
  \Delta \varphi-\omega \varphi+f(\varphi)=0 \quad\text{for $x\in\R^2$}.
\end{equation}
A standing wave solution of the form $e^{i(\omega t+m\theta)}\phi_\omega(r)$
appears in the study of nonlinear optics (see references in \cite{PWa}).
If $m=0$ and $\phi_\omega(r)$ is positive, then $\phi_\omega$ is
a ground state. Existence and uniqueness of the ground state are
well known (see \cite{BeLi1}, \cite{BLP}, \cite{Kw} and reference therein).

If $m\ne 0$, Iaia and Warchall proved the existence of smooth
solutions to \eqref{eq:B} with any prescribed number of zeroes.
The uniqueness of positive solutions has been proved by
\cite{Mi} by using the
classification theorem of positive solutions due to Yanagida and
Yotsutani \cite{YY}.

Let $c>0$ and let $Q_c$ be a positive solution to
\begin{equation}
  \label{eq:1D}
\left\{
\begin{aligned}
& Q''-cQ+f(Q)=0 \quad\text{for $x\in\R$,}
\\ &  \lim_{x\to\pm\infty}Q(x)=0,
\\ &  Q(0)=\max_{x\in\R}Q(x).
    \end{aligned}
\right.
\end{equation}
Then
\begin{equation}
  \label{eq:Qc}
Q_c(x)=\left(\frac{(p+1)c}{2}\right)^{\frac{1}{p-1}}\sech^{\frac{2}{p-1}}
\left(\frac{(p-1)\sqrt{c}}{2}x\right).
\end{equation}
\par
In \cite{PWa}, Pego and Warchall numerically observe that as spin index
$m$ becomes larger, a solution $\phi_\omega(r)$ to \eqref{eq:B} remains 
small initially and then is approximated by $Q_c(r-\bar{r})$ around 
$r=\bar{r}$, where $c=\omega+(m^2/\bar{r}^2)$
and $\bar{r}$ is a positive number with
$\bar{r}=O(m)$ as $m\to\infty$ (see also \cite{MoK} and references
 in \cite{PWa}).
One of our goals in the present paper is to explain this phenomena.
Benci and D'Aprile \cite{BdA} studied \eqref{eq:B} in a general
setting and locate the asymptotic peak of solutions (see also \cite{dA}).
Recently, Ambrosetti, Malchiodi and Ni \cite{AMN} have proved
the existence of positive radial solutions concentrating on spheres to
a class of singularly perturbed problem
$$\varepsilon^2\Delta u-Vu+|u|^{p-1}u=0,$$
and obtain their asymptotic profile.
Adopting the argument in \cite{AMN}, we obtain the following.
\begin{theorem}
\label{thm1}
Let $p>1$ and let $\phi_{\omega,m}$ be a positive solution to
\eqref{eq:B}. Then there exists an $m_*\in\N$ such that if $m\ge m_*$,
\begin{gather}
\label{eq:1.5}
\|\phi_{\omega,m}(\cdot)-Q_c(\cdot-\bar{r})\|_{H^2_{r}(\R^2)}
=O(m^{-1/2}),\\
\label{eq:1.6}
\|\phi_{\omega,m}(\cdot)-Q_c(\cdot-\bar{r})\|_{L^\infty(\R^2)}=O(m^{-1}),
\end{gather}
where $\bar{r}=2m/\sqrt{(p-1)\omega}$ and $c=(p+3)\omega/4$.
\end{theorem}
\begin{remark}
 Let $r=ms$, $\varepsilon=1/m$ and $V(r)=\omega+r^{-2}$. Then
\eqref{eq:B} is transformed into 
$$
\varepsilon^2\Delta_r \phi-V(r)\phi+f(\phi)=0.$$
Though \cite{AMN} assumes the boundedness of $V(r)$ and cannot be
applied directly to our problem,
a maximum point of $\phi_{\omega,m}(r)$ can be predicted from an
\textit{auxiliary weighted potential} $rV(r)$ introduced by \cite{AMN}. 
\end{remark}

Let $\varphi_\omega$ be a ground state to \eqref{eq:sfe}.
As is well known, the standing wave solution $e^{i\omega t}\varphi_\omega$ 
is stable if $d\|\varphi_\omega\|_{L^2(\R^n)}^2/d\omega>0$ and unstable
if $d\|\varphi_\omega\|_{L^2(\R^n)}^2/d\omega<0$. See e.g. Berestycki-Cazenave
\cite{BeCa}, Cazenave-Lions \cite{CL}, Grillakis-Shatah-Strauss
\cite{GSS1}, \cite{GSS2}, Shatah \cite{Sh},
Shatah-Strauss \cite{ShSt} and  Weinstein \cite{We2}.
Namely, the standing wave solution $e^{i\omega t}\varphi_\omega$
is stable if $1<p<1+4/n$ and unstable if $p\ge1+4/n$.
Grillakis \cite{Gr} proved that every radially symmetric standing wave
solution is linearly unstable if $p>1+4/n$.
However, to the best our knowledge, it remains unknown whether there
exists an unstable standing wave solution with higher energy 
in the subcritical case ($1<p<1+4/n$).

Using Theorem \ref{thm1}, we find an unstable direction and
prove $e^{i(\omega t+m\theta)}\phi_\omega(r)$ is unstable in $H^1(\R^2)$
if $p>1$ and $m$ is sufficiently large.
\begin{theorem}
\label{thm2}
Let $p>1$ and $\phi_{\omega,m}$ be as in Theorem \ref{thm1}.
Then there exists an $m_*\in\N$ such that if $m\ge m_*$,
a standing wave solution
$e^{i(\omega t+m\theta)}\phi_\omega$ is linearly unstable.
\end{theorem}

\begin{remark}
By Shatah-Strauss Lemma (see \cite{ShSt2,SW}, see also \cite{Mi2}),
we have orbital instability of the linearly unstable standing wave
solutions.  
\end{remark}
\begin{remark}
 If $p<1+4/n$ and $u_0\in H^1(\R^n)$, a solution to \eqref{eq:NLS}
exists globally in time and remains bounded in $H^1(\R^n)$.
Thus the mechanism of instability shown in Theorem \ref{thm2} is quite
different from that of \cite{BeCa} where solutions around a
standing wave solution blow up in finite time.
The instability mechanism we find is close to transversal long-wave
instability of 1-dimensional soliton (see Alexander-Pego-Sachs \cite{APS} for KP equation and Bridges \cite{B1,B2} for nonlinear Schr\"{o}dinger
equation). Theorem \ref{thm1} shows that a profile of vortex soliton
is close to 1D-soliton for large $m$ and thus it becomes possible
to find unstable modes by using perturbation method.
\end{remark}
\begin{remark}
If $\phi_{\omega, m}$ is nonnegative, then
$e^{i(m\theta+\omega t)}\phi_\omega(r)$ is a ground state in the class
$X_m=\{e^{im\theta}v(r)\,|\, v\in H^1_{rad}(\R^2),\;
 v\in L^2_{rad}(\R^2)\}$ and it follows from Grillakis
\textit{et al.} (\cite{GSS1}) that the standing wave solution
$e^{i(m\theta+\omega t)}\phi_\omega(r)$
is stable in the class $X_m$ if $1<p<3$ (\cite{Mi}). Thus the vortex
soliton is stable to the symmetric perturbations in the subcritical case.
\end{remark}
Our plan of the present paper is as follows.
In Section 2, we specify a solution to \eqref{eq:B}  which is expected
to become close to a solution to \eqref{eq:1D} as $m$ tends to infinity.
In Section 3, we investigate some properties of the
linearized operator around an approximate solution constructed is
Section 2. In Section 4, we prove Theorem \ref{thm1} following the
lines of \cite{AMN} and using Liapunov Schmidt method.
Since $\|\phi_\omega(r)\|_{\Lr}$ grows up as $m\to\infty$ whereas
$\|\phi_\omega(r)\|_{L^\infty}$ remains bounded,
we need to estimate both $\Lr$-norm and
$L^\infty_r(\R^2)$-norm of the solution to obtain Theorem \ref{thm1}.
In Section 5, we prove that $e^{im(\omega t+m\theta)}\phi_\omega$ is unstable
to the perturbations of the form $e^{i(m+j)\theta}v(r)$ with $|j|\sim
m^{\min(p-1,1)/6}$ and obtain Theorem \ref{thm2}.
\par
Finally, we introduce several notations.
For Banach spaces $X$ and $Y$, let $B(X,Y)$ be the space of all
bounded linear operators from $X$ to $Y$
and let $\|A\|_{B(X,Y)}$ be the operator norm of an operator
$A\colon X\to Y$. We abbreviate $B(X,X)$ as $B(X)$.
We denote by $D(A)$ and $R(A)$ the domain and the range of the
operator $A$, respectively. We use notations
$\|f\|_{\Lr}=\left(\int_0^\infty|f(r)|^2rdr\right)^{\frac12}$,
$\|f\|_{H^1_r(\R^2)}=\left(
\int_0^\infty\left(|f'(r)|^2+|f(r)|^2\right)rdr\right)^{\frac12}$,
$\Delta_r=\pd_r^2+r^{-1}\pd_r$ and
$\|f\|_{H^2_{r}(\R^2)}=\|(1-\Delta_r)f\|_{\Lr}$.
Various constants will be simply denoted by $C$ and $C_i$
($i\in\mathbb{N}$) in the course of calculations.

\section{An approximation}
\label{sec:2}
In this section, we will construct an approximate solution to
\eqref{eq:B} for large $m$.
Suppose that a positive solution to \eqref{eq:B} is approximated by
$Q_c(r-\bar{r})$ around $r=\bar{r}$ for large $m$.
Let
$\alpha_0=\bar{r}/m$, $\varepsilon=m^{-1}$, $s=r-\bar{r}$ and
$v(s)=\phi_\omega(r)$.
Then \eqref{eq:B} transforms into
\begin{equation}
  \label{eq:B1}
  \left\{
    \begin{aligned}
& v_{ss}+\frac{\varepsilon}{\alpha_0+\varepsilon s}v_s
-\left(\omega+\frac{1}{(\alpha_0+\varepsilon s)^2}\right)v+f(v)=0
\quad\text{for $s\in(-\bar{r},\infty)$,}
\\
& \lim_{s\to-\bar{r}}\frac{v(s)}{(s+\bar{r})^m}
=\lim_{s\to-\bar{r}}\frac{v_s(s)}{m(s+\bar{r})^{m-1}},
\\
& \lim_{s\to\infty}v(s)=0.
    \end{aligned}\right.
\end{equation}
Substituting $v(s)=v_0(s)+\varepsilon v_1(s)+O(\varepsilon^2)$
into \eqref{eq:B} and formally equating the power of $\varepsilon$, we 
obtain
\begin{equation}
  \label{eq:approx0}\left\{
  \begin{aligned}
  & v_0''-cv_0+f(v_0)=0,\\
  & \lim_{s\to\pm\infty}v_0(s)=0,
  \end{aligned}\right.
\end{equation}
and
\begin{equation}
  \label{eq:approx1}\left\{
  \begin{aligned}
& v_1''-cv_1+f'(v_0)v_1=-\alpha_0^{-1}v_0'-2\alpha_0^{-3}sv_0,
\\
& \lim_{s\to\pm\infty}v_1(s)=0,
  \end{aligned}\right.
\end{equation}
where $c=\omega+\alpha_0^{-2}$.
Let $v_0(s)=Q_c(s)$, $L_c:=\pd_s^2-c+f'(Q_c)$ and $D(L_c)=H^2(\R)$.
Since $\ker(L_c)=\operatorname{span}\{Q_c'\}$, the Fredholm
alternative implies that
\eqref{eq:approx1} has a solution $v_1\in L^2(\R)$ if and only if
\begin{equation}
\label{eq:alpha}
\int_\R Q_c'(s)\left(Q_c'(s)+\frac{2s}{\alpha_0^2}Q_c(s)\right)ds
=\int_\R \left(Q_c'(s)^2-\frac{1}{\alpha_0^2}Q_c(s)^2\right)ds=0.
\end{equation}
\begin{lemma}
\label{lem:amplitude}
Let $c=\omega+\alpha_0^{-2}$ and let $Q_c$ be a solution to
\eqref{eq:1D}. If \eqref{eq:alpha} holds,
then $c=(p+3)\omega/4$ and $\alpha_0=2/\sqrt{(p-1)\omega}$.
\end{lemma}
\begin{proof}
By \eqref{eq:1D},
$$
\left(\frac{dQ_c}{dx}\right)^2=cQ_c^2\left(
1-\left(\frac{Q_c}{A}\right)^{p-1}\right),
$$
where $A^{p-1}=(p+1)c/2$.
We compute
\begin{align*}
\int_{-\infty}^\infty Q_c'(x)^2dx=& 2\int_0^A 
\left(\frac{dQ_c}{dx}\right)^2\left(-\frac{dx}{dQ_c}(x)\right)dQ_c
\\ =&
2\sqrt{c}\int_0^A u\sqrt{1-\bigl(\frac{u}{A}\bigr)^{p-1}}du
\\ =&
\frac{2}{p-1}\sqrt{c}A^2B\left(\tfrac{2}{p-1},\tfrac{3}{2}\right),
\end{align*}
and
\begin{align*}
\int_{-\infty}^\infty Q_c(x)^2dx =&
2\int_0^AQ_c^2\left(-\frac{dx}{dQ_c}\right)dQ_c
\\=&
\frac{2}{\sqrt{c}}\int_0^A\frac{u}{\sqrt{1-\bigl(\frac{u}{A}\bigr)^{p-1}}}du
\\=&
\frac{2A^2}{(p-1)\sqrt{c}}B\left(\tfrac{2}{p-1},\tfrac12\right).
\end{align*}
Combining the above, we have $c=(p+3)\omega/4$  and
$\alpha_0=2/\sqrt{(p-1)\omega}$.
\end{proof}

Let $\chi(s)$ be smooth nonnegative functions on $\R$ satisfying
$0\le\chi(r)\le1$ and
$$
\chi(r)=\left\{\begin{aligned}
& 1\quad \text{ if $|r|\le 2$,}\\
& 0\quad \text{ if $|r|\ge 3$,}\end{aligned}\right.
$$
and let  $\chi_l(s)=\chi(s/l)$, where
$l=-\frac{2}{\sqrt{c}}\max(1,\frac1{p-1})\log\varepsilon$.
Following \cite{AMN}, we put
\begin{align}
\label{eq:de1}
& \Phi(\varepsilon,\rho)(r)=\chi_l(r-\rho)Q_c(r-\rho),
\quad c=\omega+(\varepsilon\rho)^{-2},
\\
\label{eq:de2}
&\phi_{\omega,m}=\Phi(\varepsilon,\rho)+w,
\end{align}
and search for a positive solution to \eqref{eq:B} for large $m$.
To fix the decomposition \eqref{eq:de2}, we assume
\begin{equation}
  \label{eq:de3}
  \left(w,\pd_\rho\Phi\right)_{\Lr}=0.
\end{equation}
Substituting \eqref{eq:de1} into \eqref{eq:B}, we obtain
\begin{equation}
  \label{eq:2.6}
\mathcal{L}(\varepsilon,\rho)w+R_1(\varepsilon,\rho,w)+R_2(\varepsilon,\rho)=0,
\end{equation}
where $R_2=R_{21}+R_{22}+R_{23}$ and
\begin{align*}
\mathcal{L}(\varepsilon,\rho)=& \Delta_r-\omega-\frac{m^2}{r^2}
+f'(\Phi(\varepsilon,\rho)),
\\
R_1=&f(\Phi(\varepsilon,\rho)+w)-f(\Phi(\varepsilon,\rho))
-f'(\Phi(\varepsilon,\rho))w,
\\
R_{21}=&f(\Phi(\varepsilon,\rho))-\tau_\rho(\chi_lf(Q_c))
\\
R_{22}=&
\left(c-\omega-\frac{m^2}{r^2}\right)\Phi(\varepsilon,\rho)
+\frac{1}{r}\tau_\rho(\chi_lQ_c')
\\
R_{23}=& \tau_\rho\left(\chi_l''Q_c+2\chi_l'Q_c'\right)
+\frac{1}{r}\tau_\rho\left(\chi_l'Q_c\right).
\end{align*}
Here $\tau_h$  denotes the translation, that is, $(\tau_hf)(x)=f(x-h)$.
We will search a solution $(\rho,w)$ to \eqref{eq:de3} and \eqref{eq:2.6}
with $\rho\in(\alpha_0/(2\varepsilon), 2\alpha_0/\varepsilon)$ for large
$m\in\N$.

\section{Spectrum of the linearized operator
$\mathcal{L}(\varepsilon,\rho)$}
In this section, we examine spectral properties of the linearized
operator $\mathcal{L}(\varepsilon,\rho)$.
To begin with, we recall some properties of the operator
$\Delta_r-\omega-m^2/r^2.$
\begin{lemma}
  \label{lem:3.1}
Let $0<\varepsilon<1/2$ and $\mathcal{L}_0(\varepsilon) \colon \Lr\to \Lr$
be a closed operator such that
$$\mathcal{L}_0(\varepsilon)u=\Delta_ru-\omega u-(\varepsilon r)^{-2}u$$
for $u\in C_0^\infty(\R_+).$
Then $\mathcal{L}_0(\varepsilon)$ is a self-adjoint operator with
$$
D(\mathcal{L}_0(\varepsilon))=\{u\in H^2_{r}(\R^2)\,|\,r^{-2}u\in \Lr\}
\quad\text{and}\quad R(\mathcal{L}_0(\varepsilon))=\Lr.$$
\end{lemma}
\begin{proof}
Let $X=\{H^2_{r}(\R^2)\,|\, r^{-2}u\in \Lr\}$ be a Hilbert space
equipped with the norm
$\|u\|_{X}=(\|u\|_{H^2_{r}(\R^2)}^2+\|r^{-2}u\|_{\Lr}^2)^{1/2}$.

By Theorem 10.10 and Example 4 in \cite[Appendix to X.1]{RS},
the operator $\mathcal{L}_0(\varepsilon)$ is essentially self-adjoint in
$C_0^\infty(\R_+)$. Thus for any $u\in D(\mathcal{L}_0(\varepsilon))$,
there exist $u_n\in C_0^\infty(\R_+)$ $(n\in\N)$ such that
$\mathcal{L}_0(\varepsilon)u_n\to \mathcal{L}_0(\varepsilon)u$ and
 $u_n\to u$ in $\Lr$ as $n\to\infty$.

Integrating by parts, we have
\begin{equation}
\label{eq:3.6}
  \begin{split}
&\|\mathcal{L}_0(\varepsilon)w\|_{\Lr}^2 \\=&
\|(\omega-\Delta_r)w\|_{\Lr}^2+
2\Re\left((\omega-\Delta_r)w,(\varepsilon r)^{-2}w\right)_{\Lr}
+\left\|(\varepsilon r)^{-2}w\right\|_{\Lr}^2
\\ \ge &
\|(\omega-\Delta_r)w\|_{\Lr}^2+
(\varepsilon^{-4}-4\varepsilon^{-2})\left\|r^{-2}w\right\|_{\Lr}^2
\end{split}
\end{equation}
for every $w\in C_0^\infty(\R_+)$.
Eq. \eqref{eq:3.6} yields that $\{u_n\}_{n=1}^\infty$ is a Cauchy
sequence in $X$ and $u_n\to u$ in $X$ as $n\to\infty$.
Thus we have $D(\mathcal{L}_0(\varepsilon))\subset X$.
\par
Next we prove $D(\mathcal{L}_0(\varepsilon))\supset X$.
For every $u\in X$, there exist $u_n\in C_0^\infty(\R_+)$ $(n=1,2,\cdots)$
such that
$\lim_{n\to\infty}\|u_n-u\|_{X}=0$.
Since
$$
\|\mathcal{L}_0(\varepsilon)w\|_{\Lr}\le \max(1,\omega)
\|w\|_{H^2_{r}(\R^2)}+\varepsilon^{-2}\|r^{-2}w\|_{\Lr},$$
we see that $\{\mathcal{L}_0(\varepsilon)u_n\}_{n=1}^\infty$
 and $\{u_n\}_{n=1}^\infty$
are Cauchy sequences in $\Lr$ and that there exists a
$v\in \Lr$ such that $\mathcal{L}_0(\varepsilon)u_n\to v$ as $n\to\infty$.
Since $\mathcal{L}_0(\varepsilon)$ is closed, it follows that
 $v=\mathcal{L}u$ and $u\in D(\mathcal{L}_0(\varepsilon))$. Thus we prove
$D(\mathcal{L}_0(\varepsilon))=X$.

Finally, we will show that $R(\mathcal{L}_0(\varepsilon))=\Lr$.
The self-adjointness of $\mathcal{L}_0(\varepsilon)$ and \eqref{eq:3.6}
implies $$R(\mathcal{L}_0(\varepsilon))^\perp
=\ker(\mathcal{L}_0(\varepsilon))=\{0\}.$$
Hence it follows that $\overline{R(\mathcal{L}_0(\varepsilon))}=\Lr$
and that for every $v\in \Lr$, there exist $u_n\in X$ and $v_n\in\Lr$
$(n\in\N)$ such that
$$\mathcal{L}_0(\varepsilon)u_n=v_n\to v\quad\text{in $\Lr$ as $n\to\infty$.}$$
By \eqref{eq:3.6}, there exists $u\in X$ such that $\lim_{n\to\infty}u_n=u$
in $X$. Since $\mathcal{L}_0(\varepsilon)$ is closed we have
$v=\mathcal{L}_0(\varepsilon)u\in R(\mathcal{L}_0(\varepsilon))$.
This completes the proof of Lemma \ref{lem:3.1}.
\end{proof}

Let $\mathcal{P}(\varepsilon,\rho)$ and $\mathcal{Q}(\varepsilon,\rho)$ be
orthogonal projections defined by
\begin{align*}
&\mathcal{P}(\varepsilon,\rho)u=\|\pd_\rho\Phi(\varepsilon,\rho)\|_{\Lr}^{-2}
(u,\pd_\rho\Phi(\varepsilon,\rho))_{\Lr}\pd_\rho\Phi(\varepsilon,\rho),
\\ & \mathcal{Q}(\varepsilon,\rho)=I-\mathcal{P}(\varepsilon,\rho).
\end{align*}
We will show that $\mathcal{L}(\varepsilon,\rho)$ is invertible on
$\mathcal{Q}(\varepsilon,\rho)\Lr$. 
\begin{lemma}
  \label{lem:3.2}
Let $w\in H^1_r(\R^2)\cap \{u\,|\, r^{-1}u\in \Lr\}$ and
$$(w,\pd_\rho\Phi(\varepsilon,\rho))_{\Lr}
=(w,\Phi(\varepsilon,\rho)^{\frac{p+1}{2}})_{\Lr}=0.$$
Then there exist positive numbers $\varepsilon_*$ and $c_1$ such that
$$-(\mathcal{L}(\varepsilon,\rho)w,w)_{\Lr} \ge c_1\|w\|_{\Lr}^2$$
for every $\varepsilon\in(0,\varepsilon_*)$ and
$\rho\in(\alpha_0/(2\varepsilon),2\alpha_0/\varepsilon)$.
\end{lemma}
To prove Lemma \ref{lem:3.2}, we need the following.
\begin{lemma}
  \label{lem:l1}
Let $p>1$, $c>0$ and $\lambda_0=(p-1)(p+3)/4$.
Then $$L_c Q_c^{\frac{p+1}{2}}=\lambda_0 cQ_c^{\frac{p+1}{2}}.$$
Furthermore, $\ker(L_c)=\{\beta Q_c'\,|\,\beta\in\R\}$
and there exists a positive number $b$ depending only on of $p$ such that
$\sigma(L_c)\setminus\{0,\lambda_0 c\} \subset (-\infty,-bc].$
\end{lemma}
\begin{proof}
The former part of the lemma can be obtained by a simple computation.
Let $c=1$. 
Weyl's essential spectrum theorem tells us that the spectrum of $L_1$ 
consists of essential spectrum $(-\infty,-1]$ and discrete eigenvalues.
Since $Q_1'$ has exactly one zero and $L_1Q_1'=0$,
it follows from Strum's comparison theorem that $0$ is a second
eigenvalue of $L_1$ and that $\ker(L_1)$ is spanned by $Q_1'$.
Since $L_c(u(c^{1/2}x))=c(L_1u)(c^{1/2}x)$ for every $u\in H^2(\R)$, we
have $ \sigma(L_c)=\{c\lambda\,|\, \lambda\in\sigma(L_1)\}.$
Thus we prove Lemma \ref{lem:l1}.
\end{proof}
\begin{proof}[Proof of Lemma \ref{lem:3.2}]
Let $\chi_0(s)=1-\chi_l(s)$ and $\chi_1(s)=\chi_l(s)$. By \eqref{eq:Qc}
and the fact that $\mathrm{supp}\chi_0\subset \{r\in\R\,|\, |r|\ge 2l\}$,
\begin{align*}
(\mathcal{L}(\varepsilon,\rho)w,w)_{\Lr}
=&
(\mathcal{L}(\varepsilon,\rho)\chi_1w,\chi_1w)_{\Lr}
+2(\mathcal{L}_0(\varepsilon)\chi_0w,\chi_1w)_{\Lr}
\\&+(\mathcal{L}_0(\varepsilon)\chi_0w,\chi_0w)_{\Lr}
+O(e^{-2(p-1)\sqrt{c}l}\|w\|_{\Lr}^2).
\end{align*}
Integrating by parts and substituting $|\chi_0'(r)|+|\chi_1'(r)|
=O(l^{-1})$ into the resulting equation, we have
\begin{equation}
  \label{eq:L00}
  \begin{split}
-&(\mathcal{L}_0(\varepsilon)\chi_0w,\chi_0w)_{\Lr}\\
=& \int_0^\infty \left((\chi_0w)_r^2+(\omega +(\varepsilon r)^{-2})(\chi_0w)^2
\right)rdr
\\=& \int_0^\infty \chi_0^2\left(w_r^2+\omega w^2+(\varepsilon r)^{-2}w^2
\right)rdr+O(l^{-1}\|w\|_{H^1_r(\R^2)}^2),
  \end{split}
\end{equation}
and
\begin{equation}
  \label{eq:L01}
  \begin{split}
&-(\mathcal{L}_0(\varepsilon)\chi_0w,\chi_1w)_{\Lr}
\\= &
\int_0^\infty \chi_0\chi_1\left(w_r^2+\omega w^2+(\varepsilon r)^{-2} w^2
\right)rdr +O(l^{-1}\|w\|_{H^1_r(\R^2)}^2).
  \end{split}
\end{equation}
Let $U\colon \Lr\to L^2(\R_+)$ be the unitary operator
defined by $U\phi(r)=r^{\frac12}\phi(r)$. 
Then
\begin{equation}
  \label{eq:L11}
  \begin{split}
\tau_{-\rho}U\mathcal{L}(\varepsilon)U^{-1}=&\pd_r^2-\omega
-\frac{1-\frac14\varepsilon^2}{(\alpha+\varepsilon r)^2}+f'(\chi_1Q_c)
\\=& L_c+\left(\frac{1}{\alpha^2}-
\frac{1-\frac14\varepsilon^2}{(\alpha+\varepsilon r)^2}\right)
+f'(\chi_1Q_c)-f'(Q_c),      
  \end{split}
\end{equation}
where $\alpha=\rho/m$ and $c=\omega+\alpha^{-2}$.
Let $\tilde{\chi}_1$ and $\tilde{\chi}_2$ be smooth nonnegative functions
on $\R$ satisfying 
\begin{align*}
& \sup_{r\in\R}|\tilde{\chi}_i'(r)|=O(l^{-1})\quad\text{for $i=0,1$,}\\  
& \tilde{\chi}_0(r)=\left\{
  \begin{aligned}
    0 &\quad\text{if $|r|\le l$}\\
    1 &\quad\text{if $|r|\ge 2l$}
  \end{aligned}\right.,\quad
\tilde{\chi}_1(r)=\left\{
  \begin{aligned}
    1 &\quad\text{if $|r|\le 3l$}\\
    0 &\quad\text{if $|r|\ge 4l$}
  \end{aligned}\right..
\end{align*}
Put $\tilde w(r)=(r+\rho)^{1/2}\chi_1(r)w(r+\rho).$
Using $w\perp \pd_\rho \Phi(\varepsilon,\rho)$ and
\begin{equation}
\label{eq:pdP}
\pd_\rho\Phi(\varepsilon,\rho)=-\tau_\rho(\chi_1Q_c)'
-\frac{2\varepsilon}{\alpha^3}\tau_\rho(\chi_1\pd_cQ_c),  
\end{equation}
we have
\begin{align*}
  0=&(w,\pd_\rho\Phi(\varepsilon,\rho))_{\Lr}\\
=& -\int_{-\rho}^\infty (\chi_1Q_c)'(r)w(r+\rho)(r+\rho)dr
+O(\varepsilon^{1/2}\|w\|_{\Lr})
\\=&
-\int_\R (\rho+r)^{1/2}\tilde{\chi}_1\tilde{w}Q_c'dr+
O\left((\rho^{1/2}e^{-2\sqrt{c}l}+\varepsilon^{1/2})\|w\|_{\Lr}\right)
\\=&-\rho^{1/2}\int_\R\tilde{w}Q_c'dr+O(\varepsilon^{1/2}l\|w\|_{\Lr}).
\end{align*}
Hence it follows that
\begin{equation}
  \label{eq:L12}
(\tilde{w},Q_c')_{L^2(\R)}=O(\varepsilon\log\varepsilon\|w\|_{\Lr}).
\end{equation}
Similarly, we have
\begin{equation}
  \label{eq:L12b}
(\tilde{w},Q_c^{\frac{p+1}{2}})_{\Lr}=O(\varepsilon\log\varepsilon\|w\|_{\Lr}).
\end{equation}
Combining Lemma \ref{lem:l1} with  \eqref{eq:L11}, \eqref{eq:L12} and
\eqref{eq:L12b}, we see that there exist positive constants $C_1$ and
$C_2$ such that
\begin{equation}
  \label{eq:L1}
-(\mathcal{L}(\varepsilon,\rho)\chi_1w,\chi_1w)_{\Lr}
\ge C_1\|\tilde w\|_{H^1(\R)}^2 \ge C_2\|\chi_1w\|_{H^1_r(\R^2)}^2.
\end{equation}
Thus by \eqref{eq:L00}, \eqref{eq:L01} and \eqref{eq:L1},
there exist positive numbers $c_1$ and $\varepsilon_*$ such that
$$-(\mathcal{L}(\varepsilon,\rho)w,w)_{\Lr}\ge c_1\|w\|_{\Lr}^2$$
for every $\varepsilon\in(0,\varepsilon_*)$,
$\rho\in(\alpha_0/(2\varepsilon), 2\alpha_0/\varepsilon)$.
\end{proof}
Let $X_1=\mathcal{Q}(\varepsilon,\rho)X$,
$Y_1=\mathcal{Q}(\varepsilon,\rho)\Lr$ and
$\mathcal{A}(\varepsilon,\rho)=
\mathcal{Q}(\varepsilon,\rho)\mathcal{L}(\varepsilon,\rho)
\mathcal{Q}(\varepsilon,\rho)$. Lemma \ref{lem:3.2} yields that
$\mathcal{A}(\varepsilon,\rho)\colon X_1\to Y_1$ is isomorphic.
\begin{corollary}
\label{cor:3.3}
There exist positive numbers  $\varepsilon_*$ and $\nu$ such that
\begin{equation}
  \label{eq:cor3.3}
\|\mathcal{A}(\varepsilon,\rho)^{-1}u\|_{X}\le \nu \|u\|_{\Lr}  
\end{equation}
for every $u\in Y_1$, $\varepsilon\in(0,\varepsilon_*)$
and $\rho\in(\alpha_0/(2\varepsilon), 2\alpha_0/\varepsilon)$.
\end{corollary}
\begin{proof}
Let $Q_1$ and $Q_2$ be orthogonal projections such that
\begin{gather*}
Q_1u=\frac{(u,\mathcal{Q}(\varepsilon,\rho)
\Phi(\varepsilon,\rho)^{\frac{p+1}{2}})_{\Lr}}
{\|\mathcal{Q}(\varepsilon,\rho)\Phi(\varepsilon,\rho)^{\frac{p+1}{2}}
\|_{\Lr}^{2}}
\mathcal{Q}(\varepsilon,\rho)\Phi(\varepsilon,\rho)^{\frac{p+1}{2}},
\\
Q_2=\mathcal{Q}(\varepsilon,\rho)-Q_1.  
\end{gather*}
Then $\mathcal{A}(\varepsilon,\rho)$ can be written as
$$
\mathcal{A}(\varepsilon,\rho)=
\begin{pmatrix}
Q_1\mathcal{L}(\varepsilon,\rho)Q_1 & Q_1\mathcal{L}(\varepsilon,\rho)Q_2\\
Q_2\mathcal{L}(\varepsilon,\rho)Q_1 & Q_2\mathcal{L}(\varepsilon,\rho)Q_2
\end{pmatrix}.$$
In view of Lemma \ref{lem:l1}, we see that there exists a $c_2>0$ such that
$$
(\mathcal{L}(\varepsilon,\rho)\mathcal{Q}(\varepsilon,\rho)
\Phi(\varepsilon,\rho)^{\frac{p+1}{2}},\mathcal{Q}(\varepsilon,\rho)
\Phi(\varepsilon,\rho)^{\frac{p+1}{2}})_{\Lr}
\ge c_2
\|\mathcal{Q}(\varepsilon,\rho)\Phi(\varepsilon,\rho)^{\frac{p+1}{2}}
\|_{\Lr}^2.$$
Furthermore, we see that
$$
\lim_{\varepsilon\downarrow0}
\left(\|Q_1\mathcal{L}(\varepsilon,\rho)Q_2\|_{B(\Lr)}
+\|Q_2\mathcal{L}(\varepsilon,\rho)Q_1\|_{B(\Lr)}\right)=0.$$
Combining the above with Lemma \ref{lem:3.2}, we obtain
\begin{equation}
\label{eq:3.7}
\sup_{\varepsilon\in(0,\varepsilon_*)}\sup_{\rho\in
(\alpha_0/(2\varepsilon), 2\alpha_0/\varepsilon)}
\|\mathcal{A}(\varepsilon,\rho)^{-1}\|_{B(\Lr)}<\infty.
\end{equation}

Let 
\begin{align*}
\mathcal{B}(\varepsilon,\rho)=&
\mathcal{P}(\varepsilon,\rho)\mathcal{L}(\varepsilon,\rho)
+\mathcal{L}(\varepsilon,\rho)\mathcal{P}(\varepsilon,\rho)
-\mathcal{P}(\varepsilon,\rho)\mathcal{L}(\varepsilon,\rho)
\mathcal{P}(\varepsilon,\rho) -f'(\Phi(\varepsilon,\rho)).
\end{align*} Then
\begin{equation}
  \label{eq:3.8}
\mathcal{L}_0(\varepsilon)=\mathcal{A}(\varepsilon,\rho)
+\mathcal{B}(\varepsilon,\rho).
\end{equation}
Using \eqref{eq:3.7}, \eqref{eq:3.8} and the fact that
$$
\sup_{\varepsilon\in(0,\varepsilon_*)}\sup_{\rho\in
(\alpha_0/(2\varepsilon), 2\alpha_0/\varepsilon)}
\left\|\mathcal{B}(\varepsilon,\rho)\right\|_{B(\Lr)}<\infty,$$
we have
\begin{equation}
  \label{eq:3.9}
\begin{split}
\|\mathcal{L}_0(\varepsilon)\mathcal{A}(\varepsilon,\rho)^{-1}u\|_{\Lr}
\le C\|u\|_{\Lr}
\end{split}  
\end{equation}
for every $u\in Y_1$, $\varepsilon\in(0,\varepsilon_*)$
and $\rho\in(\alpha_0/(2\varepsilon), 2\alpha_0/\varepsilon)$.
Combining \eqref{eq:3.6} and \eqref{eq:3.9}, we obtain \eqref{eq:cor3.3}.
\end{proof}
We will use the lemma below to estimate $L^\infty$-norm of $w$
in the following section.
\begin{corollary}
  \label{cor:3.4}
Let $p>1$. Then there exist positive numbers $\varepsilon_*$ and $C$
such that
\begin{equation}
\label{eq:cor34}
 \|\mathcal{A}(\varepsilon,\rho)^{-1}u\|_{L_r^\infty(\R^2)}
\le C\|u\|_{L_r^\infty(\R^2)}
\end{equation}
for every $u\in L^\infty_r(\R^2)\cap Y_1$,
$\varepsilon\in(0,\varepsilon_*)$ with $\varepsilon^{-1}\in \N$ and $\rho\in
(\alpha_0/(2\varepsilon),2\alpha_0/\varepsilon)$.
\end{corollary}
\begin{proof}
Let $m=\varepsilon^{-1}\in\N$ and
\begin{align*}
& P^{\perp}u=u-\|Q_c'\|_{L^2(\R)}^{-2}(u,Q_c')_{L^2(\R)}Q_c',\\
& \mathcal{K}(\varepsilon,\rho)
=\mathcal{Q}(\varepsilon,\rho)\left\{
(\tau_\rho\tilde{\chi}_0)\mathcal{L}_0(\varepsilon)^{-1}(\tau_\rho\chi_0)
+U^{-1}\tau_\rho\tilde{\chi}_1P^\perp L_c^{-1}P^\perp\chi_1\tau_{-\rho}
U\right\}\mathcal{Q}(\varepsilon,\rho).
\end{align*}
Noting that $e^{im\theta}\mathcal{L}_0(\varepsilon)u(r)
=(\Delta-\omega)(e^{im\theta}u(r))$, we have
$$\sup_{m\in\N}\|\mathcal{L}_0(\varepsilon)^{-1}\|_{B(L^\infty_r(\R^2))}
<\infty.$$
Furthermore, $L_c\colon P^\perp L^2(\R)\to L^2(\R)$ has a bounded
 inverse. Hence it follows that
\begin{equation}
  \label{eq:d1}
\sup_{m\ge\varepsilon_*^{-1}}\sup_{\rho\in
(\alpha_0/(2\varepsilon), 2\alpha_0/\varepsilon)}
\|\mathcal{K}(\varepsilon,\rho)\|_{B(L^\infty_r(\R^2))}<\infty.
\end{equation}
We compute
\begin{equation}
  \label{eq:d2}
\begin{split}
& \mathcal{A}(\varepsilon,\rho)\mathcal{K}(\varepsilon,\rho)
\\ =&
\mathcal{A}(\varepsilon,\rho)(\tau_\rho\tilde{\chi}_0)
\mathcal{L}_0(\varepsilon)^{-1}(\tau_\rho\chi_0)
+\mathcal{A}(\varepsilon,\rho) U^{-1}\tau_\rho
\tilde{\chi}_1P^\perp L_c^{-1}P^\perp\chi_1\tau_{-\rho} U.
\\=&I+II.
\end{split}  
\end{equation}
Since $\mathcal{A}(\varepsilon,\rho)=\mathcal{L}_0(\varepsilon)-
\mathcal{B}(\varepsilon,\rho)$ and
$\|\mathcal{B}(\varepsilon,\rho)\tau_\rho\tilde{\chi}_0\|_{B(L^\infty)}=
O(e^{-\tilde{p}\sqrt{c}l})$, where $\tilde{p}=\min(1,p-1)$, we have
\begin{equation}
  \label{eq:d3}
\begin{split}
I=& \mathcal{Q}(\varepsilon,\rho)\mathcal{L}_0(\varepsilon)
(\tau_\rho\tilde{\chi}_0)\mathcal{L}_0(\varepsilon)^{-1}(\tau_\rho\chi_0)
-\mathcal{B}(\varepsilon,\rho)
(\tau_\rho\tilde{\chi}_0)\mathcal{L}_0(\varepsilon)^{-1}(\tau_\rho\chi_0)
\\=&
\mathcal{Q}(\varepsilon,\rho)\left\{\tau_\rho(\tilde{\chi}_0\chi_0)
+[\Delta_r,\tau_\rho\tilde{\chi}_0]\mathcal{L}_0(\varepsilon)^{-1}
\tau_\rho\chi_0
-\mathcal{B}(\varepsilon,\rho)(\tau_\rho\tilde{\chi}_0)
\mathcal{L}_0(\varepsilon)^{-1}(\tau_\rho\chi_0)\right\}
\\=&\tau_\rho\chi_0+O(l^{-1})\quad\text{in
 $B(\mathcal{Q}(\varepsilon,\rho)L^\infty_r(\R^2))$.}
\end{split}  
\end{equation}
Let
$\mathcal{B}_1=
\mathcal{P}(\varepsilon,\rho)
\mathcal{L}(\varepsilon,\rho)\mathcal{P}(\varepsilon,\rho)
-\mathcal{L}(\varepsilon,\rho)\mathcal{P}(\varepsilon,\rho).$
Then $$\mathcal{A}(\varepsilon,\rho)=
\mathcal{Q}(\varepsilon,\rho)\mathcal{L}(\varepsilon,\rho)
+\mathcal{B}_1(\varepsilon,\rho).$$
In view of the definition of $\mathcal{P}(\varepsilon,\rho)$, \eqref{eq:pdP}
and the fact that $L_cQ_c'=0$, we have 
$$\left\|\mathcal{B}_1(\varepsilon,\rho)\right\|_{B(L_r^\infty(\R^2))}
=O(\varepsilon l)$$
for $\varepsilon\in(0,\varepsilon_*)$ and $\rho\in
(\alpha_0/(2\varepsilon),2\alpha_0/\varepsilon)$.
Furthermore \eqref{eq:pdP} implies
$$
\left\|\mathcal{Q}(\varepsilon,\rho)-P^\perp\right\|_{B(L_r^\infty(\R^2))}
=O(\varepsilon).$$
Let $$\mathcal{R}=\alpha^{-2}-(\alpha+\varepsilon r)^{-2}
+f'(\chi_1Q_c)-f'(Q_c).$$
Then we have $\mathcal{L}(\varepsilon,\rho)=U^{-1}(\tau_\rho L_c)U
+\tau_\rho\mathcal{R}$ and
$$
\|(\tau_\rho\mathcal{R}\tilde{\chi}_1)(1-\Delta_r)^{-1}
\|_{B(L^\infty(\R^2))}
=O(\varepsilon l+e^{-2(p-1)\sqrt{c}l}).$$
Combining the above, we have
\begin{equation}
\label{eq:d4}
\begin{split}
II=& \mathcal{Q}(\varepsilon,\rho)U^{-1} \tau_\rho(P^\perp L_c)
\tau_\rho\tilde{\chi}_1 P^\perp L_c^{-1}P^\perp\chi_1\tau_{-\rho}U
+O(\varepsilon l)
\\=&
\mathcal{Q}(\varepsilon,\rho)\tau_\rho(\tilde{\chi}_1\chi_1)+O(l^{-1})
\\=& \tau_\rho\chi_1+O(l^{-1})\quad
\text{in $B(\mathcal{Q}(\varepsilon,\rho)L^\infty_r(\R^2)$.}
\end{split}
\end{equation}
From \eqref{eq:d1}--\eqref{eq:d4}, we deduce
\eqref{eq:cor34}.
\end{proof}
\section{The method of Liapunov-Schmidt}
In this section, we use the method of Liapunov-Schmidt to obtain a
solution to \eqref{eq:de3} and \eqref{eq:2.6}.
Let us translate \eqref{eq:2.6} into a system
\begin{align}
& \label{eq:4.1}
\mathcal{A}(\varepsilon,\rho)w+\mathcal{Q}(\varepsilon,\rho)
R_1(w,\varepsilon,\rho)+\mathcal{Q}(\varepsilon,\rho)R_2(\varepsilon,\rho)=0,
\\ &
\label{eq:4.2}
\mathcal{P}(\varepsilon,\rho)\left(\mathcal{L}(\varepsilon,\rho)w
+R_1(w,\varepsilon,\rho)+R_2(\varepsilon,\rho)\right)=0.
\end{align}
\begin{lemma}
  \label{lem:4.1}
Let $p>1$. Then there exist an $\varepsilon_0>0$ and a $C>0$ such that if 
$\varepsilon\in(0,\varepsilon_0]$ and $\rho\in(\alpha_0/(2\varepsilon),
2\alpha_0/\varepsilon)$,
Eq. \eqref{eq:4.1} has a unique solution
$w(\varepsilon,\rho)$ that is continuous in $\varepsilon$ and $\rho$
and satisfies
\begin{gather}
  \label{eq:4.3}
\|w(\varepsilon,\rho)\|_X\le C\varepsilon^{\frac12}\quad\text{as 
$\varepsilon\downarrow0$.}
\end{gather}
\end{lemma}
\begin{proof}
Let $T\colon X_1\times(0,\varepsilon_0]\times(\alpha_0/(2\varepsilon),
2\alpha_0/\varepsilon)\to X_1$ be a continuous mapping defined by
$$T(w,\varepsilon,\rho)=-\mathcal{A}(\varepsilon,\rho)^{-1}
\mathcal{Q}(\varepsilon,\rho)
\left\{R_1(w,\varepsilon,\rho)+R_2(\varepsilon,\rho)\right\},
$$ and let
$\widetilde{X}=\{w\in X_1\,|\,\|w\|_{X}\le r_0\},$
where $r_0$ is a positive number to be fixed later.

To begin with, we will show that $T$ maps $\widetilde{X}$ into itself.
We compute
\begin{equation}
\label{eq:4.4}
\begin{split}
\|R_{1}\|_{\Lr} =& \left\|\int_0^1\left\{f'(\Phi(\varepsilon,\rho)+\theta w)
-f'((\Phi(\varepsilon,\rho))\right\}d\theta w\right\|_{\Lr}
\\ \le & \delta(r_0)\|w\|_{\Lr},
\end{split}
\end{equation}
where $\delta(r_0)$ is a positive constant with
$\lim_{r_0\downarrow0}\delta(r_0)=0$.
Eq. \eqref{eq:Qc} and the definition of $\chi_l$ imply
\begin{equation*}
  \begin{split}
\|R_{21}\|_{\Lr} =& \left\|\tau_\rho\left\{(\chi_1^{p-1}-\chi_1)
Q_c^{p-1}\right\}\right\|_{\Lr}
\\ \le & C\rho^{1/2}e^{-2(p-1)\sqrt{c}l},
  \end{split}
\end{equation*}
and
$$ \|R_{23}\|_{\Lr} \le C\rho^{1/2}e^{-2\sqrt{c}l}.$$
Since $\rho^{-1}=O(\varepsilon)$ and
$l=-\frac{2}{\sqrt{c}}\max(1,\frac1{p-1})\log\varepsilon$,
\begin{equation}
  \label{eq:4.5}
\|R_{21}\|_{\Lr}+ \|R_{23}\|_{\Lr}\le C_1\varepsilon^{\frac72}.
\end{equation}
Using \eqref{eq:Qc} and $\alpha^{-2}-(\alpha+\varepsilon
s)^{-2}=\frac{2s}{\alpha^3}\varepsilon+O(\varepsilon^2 s^2),$
we have
\begin{align*}
& \left\|\left(c-\omega-\frac{m^2}{r^2}\right)\tau_\rho(\chi_1Q_c')
\right\|_{\Lr}^2
\\ = & \varepsilon^{-1}\int_{-3l}^\infty(\alpha+\varepsilon s)
\left\{\left(\frac{1}{\alpha^2}-\frac{1}{(\alpha+\varepsilon s)^2}\right)
\chi_1(s)Q_c'(s)\right\}^2ds
\\ \le & C\varepsilon
\end{align*}
for every $\alpha=\varepsilon\rho\in(\alpha_0/2,2\alpha_0)$.
Similarly, we have
$$\left\|\frac{1}{r}\chi_1Q_c'\right\|_{\Lr}\le C\varepsilon^{1/2}.$$
Thus we obtain
\begin{equation}
  \label{eq:4.6}
\|R_{22}\|_{\Lr}\le C_2\varepsilon^{1/2}.
\end{equation}
Combining \eqref{eq:4.4}--\eqref{eq:4.6} with Corollary \ref{cor:3.3}, we have
\begin{equation}
  \label{eq:4.7}
 \|T(w,\varepsilon,\rho)\|_{X}\le
\nu(\delta(r_0)\|w\|_{X}+C_1\varepsilon^{7/2}+C_2\varepsilon^{1/2}).
\end{equation}
Put $r_0=2\nu C_2\varepsilon^{1/2}$. Then $T(\cdot,\varepsilon,\rho)$
maps $\widetilde{X}$ into itself if  $\varepsilon_0$ is  sufficiently
small.
\par
Next, we will show that $T(\cdot,\varepsilon,\rho)$ is a contraction mapping.
For $w_1$, $w_2\in \widetilde{X}$, 
\begin{align*}
& \|T(w_1,\varepsilon,\rho)-T(w_2,\varepsilon,\rho)\|_{X}
\\ \le & 
\nu\|R_1(w_1,\varepsilon,\rho)-R_1(w_2,\varepsilon,\rho)\|_{\Lr}
\\= &
\nu\left\|\int_0^1\left\{f'(\Phi+\theta w_1+(1-\theta)w_2)-f'(\Phi)\right\}
d\theta(w_1-w_2)\right\|_{\Lr}
\\ \le &
\Lambda\|w_2-w_1\|_{\Lr}
\end{align*}
where
$\Lambda=\nu r_0^{\tilde{p}}
\sup_{\eta\in\widetilde{X}}\|f'(\Phi+\eta)\|_{C^{\tilde{p}}}$
and $\tilde{p}=\min(1,p-1)$.
Taking $\varepsilon_0$ smaller if necessary,
we see that $T(\cdot,\varepsilon,\rho)\colon
\widetilde{X}\to\widetilde{X}$ is a contraction mapping.
Thus we prove that there exists a solution $w(\varepsilon,\rho)$ to
 \eqref{eq:4.1} with $\|w\|_{X}\le 2\nu C_2\varepsilon^{1/2}$
that is continuous in $\varepsilon\in(0,\varepsilon_0)$ and $\rho\in
(\alpha_0/(2\varepsilon),2\alpha_0/\varepsilon)$.
\end{proof}

\begin{corollary}
  \label{cor:4.2}
Let $p>1$. Then there exist an $\varepsilon_0>0$ and a $C>0$ such that
if $\varepsilon\in(0,\varepsilon_0)$, $\varepsilon^{-1}\in\N$ and
$\rho\in (\alpha_0/(2\varepsilon),2\alpha_0/\varepsilon)$,
a solution $w(\varepsilon,\rho)$ to \eqref{eq:4.1} satisfies 
\begin{equation}
  \label{eq:4.8}
\|w(\varepsilon,\rho)\|_{L^\infty} \le C\varepsilon.  
\end{equation}
\end{corollary}
\begin{proof}
Analogously to \eqref{eq:4.4}--\eqref{eq:4.6}, we have
\begin{align*}
&\|R_{1}\|_{L^\infty}\le \delta(r_0)\|w\|_{L^\infty},
\\
&\|R_2\|_{\Lr}\le \|R_{21}\|_{L^\infty}+\|R_{22}\|_{L^\infty}
+\|R_{23}\|_{L^\infty}=O(\varepsilon),
\end{align*}
where $\delta(r_0)$ is a positive number with
$\lim_{r_0\downarrow0}\delta(r_0)=0$.
Thus by Corollary \ref{cor:3.4},
\begin{align*}
 \|w(\varepsilon,\rho)\|_{L^\infty} \le & C\left(\|R_1\|_{L^\infty}+
\|R_2\|_{L^\infty}\right)
\\ \le &
C\delta(r_0)\|w\|_{L^\infty}+O(\varepsilon).
\end{align*}
Thus we have \eqref{eq:4.8}
\end{proof}
Let
$$
F(\varepsilon,\rho)=\left(\mathcal{L}(\varepsilon,\rho)w(\varepsilon,\rho)
+R_1(w(\varepsilon,\rho),\varepsilon,\rho)+R_2(\varepsilon,\rho),
\pd_\rho\Phi(\varepsilon,\rho)\right)_{\Lr}.$$
By Lemma \ref{lem:4.1}, the system of \eqref{eq:4.1} and \eqref{eq:4.2}
is reduced to an equation
\begin{equation}
   \label{eq:4.10}
F(\varepsilon,\rho)=0.
\end{equation}
\begin{lemma}
  \label{lem:4.2}
Let $p>1$ and let $\varepsilon_0>0$ be a sufficiently small number.
If $\varepsilon\in(0,\varepsilon_0]$, there exists a $\rho=\rho(\varepsilon)\in
(\alpha_0/(2\varepsilon),2\alpha_0/\varepsilon)$ satisfying
\eqref{eq:4.10}.
\end{lemma}
\begin{proof}
Let
$\mathcal{R}_c=
\frac{\varepsilon}{\alpha+\varepsilon r}\pd_r+
\left(\frac{1}{\alpha^2}-\frac{1}{(\alpha+\varepsilon r)^2}\right)
+f'(\chi_1Q_c)-f'(Q_c).$
Using \eqref{eq:Qc}, the definition of $\chi_1$ and the fact that 
$L_cQ_c'=0$ and $\rho=O(\varepsilon^{-1})$, we compute
\begin{equation}
  \label{eq:4.11}
  \begin{split}
& \|\mathcal{L}(\varepsilon,\rho)\tau_\rho(\chi_1Q_c')\|_{\Lr}
\\ \le & \|\tau_\rho\chi_1L_cQ_c'\|_{\Lr}
+\|[\pd_r^2,\tau_\rho\chi_1]\tau_\rho Q_c'\|_{\Lr}
+\left\|\tau_\rho(\mathcal{R}_c(\chi_1Q_c'))\right\|_{\Lr}
\\ = & O(\varepsilon^{1/2}).
\end{split}
\end{equation}
Similarly, we have
\begin{equation}
  \label{eq:4.12}
  \begin{split}
& \|\mathcal{L}(\varepsilon,\rho)(\pd_\rho\Phi(\varepsilon,\rho)+
\tau_\rho(\chi_1Q_c'))\|_{\Lr}\\=&
 \left\|\mathcal{L}(\varepsilon,\rho)\tau_\rho\left(\chi_1'Q_c+
\frac{2}{\alpha^3}\varepsilon\chi_1\pd_cQ_c\right)\right\|_{\Lr}
\\ = &O(\varepsilon^{1/2}).
  \end{split}
\end{equation}
By Lemma \ref{lem:4.1}, \eqref{eq:4.11} and \eqref{eq:4.12},
\begin{equation}
  \label{eq:4.13}
  \begin{split}
\left|(\mathcal{L}(\varepsilon,\rho)w,\pd_\rho\Phi(\varepsilon,\rho))_{\Lr}
\right| \le & C\varepsilon^{1/2}\|w\|_{\Lr}=O(\varepsilon).
  \end{split}
\end{equation}
Lemma \ref{lem:4.1} and Corollary \ref{cor:4.2} yield
\begin{equation}
  \label{eq:4.14}
\|R_1\|_{\Lr}\le C\|w\|_{L^\infty}^{\tilde{p}}\|w\|_{\Lr}
=O(\varepsilon^{\tilde{p}+\frac12}),
\end{equation}
where $\tilde p=\min(p-1,1)$.
Combining \eqref{eq:4.5} and \eqref{eq:4.14} with
$$\|\pd_\rho\Phi(\varepsilon,\rho)\|_{\Lr}=O(\varepsilon^{-1/2}),$$ we have
\begin{equation}
  \label{eq:4.15}
\left|(R_1+R_{21}+R_{23},\pd_\rho\Phi)_{\Lr}\right|=O(\varepsilon^{\tilde{p}}).
\end{equation}
In view of \eqref{eq:4.6} and the fact that
$\|\pd_\rho\Phi(\varepsilon,\rho)+\tau_\rho(\chi_1Q_c')\|_{\Lr}
=O(\varepsilon^{1/2}),$
\begin{equation}
  \label{eq:4.16}
 (R_{22},\pd_\rho\Phi(\varepsilon,\rho)+\tau_\rho(\chi_1Q_c'))_{\Lr}  
=O(\varepsilon).
\end{equation}
By \eqref{eq:4.13}, \eqref{eq:4.15} and \eqref{eq:4.16},
$$F(\varepsilon,\rho)=-(R_{22},\tau_\rho(\chi_1Q_c'))_{\Lr}
+O(\varepsilon^{\tilde{p}}).$$
Substituting
$$\frac{1}{\alpha^2}-\frac1{(\alpha+\varepsilon s)^2}
=\frac{2\varepsilon}{\alpha^3}s +O(\varepsilon^2s^2) 
\quad\text{as $\varepsilon\downarrow0$,}$$
and integrating by parts, we have
\begin{align*}
& (R_{22},\tau_\rho(\chi_1Q_c'))_{\Lr}
\\=&
\frac{1}{\varepsilon}\int_{-\rho}^\infty \left(\frac{1}{\alpha^2}-
\frac{1}{(\alpha+\varepsilon s)^2}\right)\chi_1(s)^2Q_c(s)Q_c'(s)
(\alpha+\varepsilon s)ds\\
& +\int_{-\rho}^\infty \chi_1(s)^2Q_c'(s)^2ds
\\ =&
\int_{-\rho}^\infty\chi_1^2
\left(\frac{2s}{\alpha^2}Q_cQ_c'+Q_c'^2\right)ds+O(\varepsilon)
\\=&
\int_\R \left\{Q_c'^2-\frac1{\alpha^2}Q_c^2\right\}ds+O(\varepsilon).
\end{align*}
Combining the above,  we see that 
$$
F(\varepsilon,\rho)=\int_\R
\left(Q_c'(s)^2-(\varepsilon\rho)^{-2}Q_c(s)^2\right)ds
+O(\varepsilon^{\tilde{p}}),
$$
where $c=\omega+(\varepsilon\rho)^{-2}$.
Hence it follows from Lemma \ref{lem:amplitude} and the intermediate
value theorem that \eqref{eq:4.10} has a solution
$\rho=\rho(\varepsilon)$ satisfying
$$\rho=(\alpha_0+o(1))\varepsilon^{-1} \quad\text{as $\varepsilon\downarrow0$.}
$$
Thus we complete the proof of Lemma \ref{lem:4.2}.
\end{proof}
Now, we are in position to prove Theorem \ref{thm1}.
\begin{proof}[Proof of Theorem \ref{thm1}]
Lemmas \ref{lem:4.1} and \ref{lem:4.2} and Corollary \ref{cor:4.2} imply
that there exists a solution $\phi_\omega$ to \eqref{eq:B} satisfying
\eqref{eq:1.5} and \eqref{eq:1.6}.
Suppose that $\phi_\omega$ is a sign-changing solution.
Since $\phi_\omega''\ge0$ and $\phi_\omega'=0$ at the minimum point,
it follows from \eqref{eq:B} that
$$\min_{r>0}\phi_\omega(r)<-\omega^{1/(p-1)}.$$
But this contracts to \eqref{eq:1.6} if $\varepsilon>0$ is
sufficiently small. Thus the solution $\phi_\omega$ to \eqref{eq:B} is
nonnegative. Since a nonnegative solution is unique (see \cite{Mi}),
we obtain Theorem \ref{thm1}.
\end{proof}
\bigskip

\section{Instability of vortex solitons}
In this section, we will prove Theorem \ref{thm2}.
Let $u(x,t)=e^{i\omega t}(e^{im\theta}\phi_\omega(r)+e^{\lambda t}v)$ and
linearize \eqref{eq:NLS} around $v=0$ and $t=0$. Then
\begin{equation}
\label{eq:5.1}
i\lambda v+(\Delta-\omega+\beta_1(r))v
+e^{2im \theta}\beta_2(r)\bar{v}=0,
\end{equation}
where
$$\beta_1(r)=\frac{p+1}2\phi_\omega(r)^{p-1},\quad
\beta_2(r)=\frac{p-1}2\phi_\omega(r)^{p-1}.$$
Put $v=e^{i(j+m)\theta}y_+$, $\bar{v}=e^{i(j-m)\theta}y_-$ and
complexify \eqref{eq:5.1} into a system
\begin{equation}
  \label{eq:5.2}
\left\{
\begin{aligned}
&\left(\Delta_r-\omega-\frac{(m+j)^2}{r^2}+i\lambda+\beta_1(r)
\right)y_+ +\beta_2(r)y_-=0,
\\
& \left(\Delta_r-\omega-\frac{(m-j)^2}{r^2}-i\lambda+\beta_1(r)
\right)y_- +\beta_2(r)y_+=0.
\end{aligned}\right.
\end{equation}
If $\lambda$ is an eigenvalue of the linearized operator,
there exist a $j\in\Z$ and a solution $(y_+,y_-)$ to
\eqref{eq:5.2} that satisfy
$(e^{i(j+m)\theta}y_+(r), e^{i(j-m)\theta}y_-(r)) \in H^1(\R^2,\C^2)$.
We will show the existence of unstable eigenvalues for $j$ with $1\ll j\ll m$.

Let $w_1=y_++y_-$, $w_2=y_+-y_-$, $\varepsilon=m^{-1}$ and
$\delta=j\varepsilon$. Let $s=r-\alpha_0m$.
Then \eqref{eq:5.2} can be rewritten as
\begin{equation}
   \label{eq:5.3}
\mathcal{H}(\varepsilon,\delta)\mathbf{w}=\lambda \mathbf{w},
\end{equation}
where $\mathbf{w}={}^t\!(w_1,w_2)$,
\begin{equation}
\label{eq:5.4}
\mathcal{H}(\varepsilon,\delta)=i
\begin{pmatrix}
 h_{11} & h_{12}\\
h_{21} & h_{22}
\end{pmatrix},
\end{equation}
and
\begin{align*}
h_{11}&= h_{22}=\frac{-2mj}{r^2},
\\
h_{12}&=\Delta_r-\omega-\frac{m^2+j^2}{r^2}+\phi_\omega^{p-1}
\\
h_{21}&=\Delta_r-\omega-\frac{m^2+j^2}{r^2}+p\phi_\omega^{p-1}.
\end{align*}
We remark that
\begin{align*}
\tau_{-\bar{r}}h_{11}=&\tau_{-\bar{r}}h_{22}=
\frac{-2\delta}{(\alpha_0+\varepsilon r)^2}
\\ 
\tau_{-\bar{r}}h_{12}=&
\pd_r^2+\frac{\varepsilon}{\alpha_0+\varepsilon r}\pd_r
-\omega-\frac{1+\delta^2}{(\alpha_0+\varepsilon r)^2}+\phi_\omega^{p-1}
\\ 
\tau_{-\bar{r}}h_{21}=&
\pd_r^2+\frac{\varepsilon}{\alpha_0+\varepsilon r}\pd_r
-\omega-\frac{1+\delta^2}{(\alpha_0+\varepsilon r)^2}+p\phi_\omega^{p-1}.
\end{align*}

Before we investigate the spectrum of $\mathcal{H}(\varepsilon,\delta)$,
let us consider the spectrum of a linear operator
$$H(\delta):=i
\begin{pmatrix}
  -2\alpha_0^{-2}\delta & L_--\alpha_0^{-2}\delta^2
\\ L_+-\alpha_0^{-2}\delta^2 & -2\alpha_0^{-2}\delta,
\end{pmatrix}$$
where $L_+=\pd_s^2-c+pQ_c^{p-1}$,  $L_-=\pd_s^2-c+Q_c^{p-1}$,
$D(L_+)=D(L_-)=H^2(\R)$ and $c=\omega+\alpha_0^{-2}$.

To begin with, we recall some spectral properties of $H(0)$.
Let
$$
\Phi_1=\begin{pmatrix}0 \\Q_c\end{pmatrix}, \quad
\Phi_2=-i\begin{pmatrix}\pd_cQ_c \\0\end{pmatrix}, \quad
\Phi_3=\begin{pmatrix}Q_c'\\0\end{pmatrix},\quad
\Phi_4=-\frac{i}{2}\begin{pmatrix}0\\sQ_c\end{pmatrix},
$$
and
$$
\Phi_1^*=\theta_1\sigma_2\Phi_2,\quad
\Phi_2^*=\theta_1\sigma_2\Phi_1,\quad
\Phi_3^*=\theta_2\sigma_2\Phi_4,\quad
\Phi_4^*=\theta_2\sigma_2\Phi_3,
$$
where
$$
\sigma_2=\begin{pmatrix}  0 & -i \\ i & 0\end{pmatrix},\quad
\theta_1=2\left(\frac{d}{dc}\|Q_c\|_{L^2(\R)}^2\right)^{-1},
\quad \theta_2=4\|Q_c\|_{L^2(\R)}^{-2}.
$$
Then we have
\begin{gather}
\label{eq:5.5}
H(0)\Phi_1=0, \quad H(0)\Phi_2=\Phi_1, \quad H(0)\Phi_3=0,
\quad H(0)\Phi_4=\Phi_3,
\\
H(0)^*\Phi_1^*=\Phi_2^*, \quad H(0)^*\Phi_2^*=0,\quad
  H(0)^*\Phi_3^*=\Phi_4^*, \quad H(0)^*\Phi_4^*=0,
\end{gather}
and $\langle \Phi_i, \Phi_j^*\rangle=\delta_{ij}$ for $i$, $j=1,2,3,4$.
Here we denote by $\langle\cdot,\cdot\rangle$ the inner product of
$L^2(\R,\C^2)$.
\begin{proposition}[see \cite{We1}]
  \label{prop:5.1}
Let $p>1$ and $p\ne 5$. Then $\lambda=0$ is a discrete eigenvalue of
$H(0)$ with algebraic multiplicity $4$.
\end{proposition}

Using Proposition \ref{prop:5.1}, we investigate the spectrum of $H(\delta)$.
\begin{lemma}
  \label{lem:5.2}
Let $1<p<5$. Then there exist a positive number $\delta_0$ and a neighborhood
$U\subset\C$ of $0$ such that for every $\delta\in(0,\delta_0)$,
$\sigma(H(\delta))\cap U$ consists of algebraically simple eigenvalues
$\lambda_i(\delta)$ $(i=1,2,3,4)$ satisfying
$$
\left|\Re\lambda_1(\delta)-\alpha_0^{-1}\gamma\delta\right|
\le \alpha_0^{-1}\gamma\delta/4,
\quad
\liminf_{\delta\downarrow0}\left(
\delta^{-1}\min_{\substack{1\le i,j\le 4,\\ i\ne j}}
|\lambda_i(\delta)-\lambda_j(\delta)|\right)
>0,$$ where
$$
\gamma=\left(2\frac{\|Q_c\|_{L^2(\R)}^2}{{\frac{d}{dc}\|Q_c\|_{L^2(\R)}^2}}
\right)^{1/2}.$$
\end{lemma}
\begin{proof}
Let $P_{H}(\delta)$ be a projection defined by
$$
P_H(\delta)=\frac{1}{2\pi i}\oint_{|\lambda|=\rho_0}
\left(\lambda-H(\delta)\right)^{-1}d\lambda,
$$
and let $Q_H(\delta)=I-P_H(\delta)$.
In view of Proposition \ref{prop:5.1}, there exist positive numbers
$\rho_0$ and $\delta_0$ such that $\mathcal{X}_0:=R(P_H(\delta))$ is
4-dimensional for every $\delta\in(0,\delta_0)$.
\par
Let $\mathcal{X}_0$ be a linear subspace whose basis is 
$\la\Phi_1,\Phi_2,\Phi_3,\Phi_4\ra$.
We decompose $H^2(\R;\C^2)$ and $L^2(\R;\C^2)$ as
$$
H^2(\R;\C^2)=\mathcal{X}_0\oplus Q_H(0)H^2(\R;\C^2),\quad
L^2(\R;\C^2)=\mathcal{X}_0\oplus Q_H(0)L^2(\R;\C^2).
$$
Then  $$H(\delta) =
\begin{pmatrix}
H_{11}(\delta) & H_{12}(\delta) \\
H_{21}(\delta) & H_{22}(\delta)
\end{pmatrix},$$
where
\begin{align*}
&H_{11}(\delta)=P_H(0)H(\delta)P_H(0),\quad
H_{12}(\delta)=P_H(0)H(\delta)Q_H(0)\\
&H_{21}(\delta)=Q_H(0)H(\delta)P_H(0),\quad
H_{22}(\delta)=Q_H(0)H(\delta)Q_H(0).
\end{align*}
By a simple computation, we have
$$H_{11}(\delta)=-2i\alpha_0^{-2}\delta I+
\begin{pmatrix}
0 & 1+b_2\delta^2 & 0 & 0\\
b_1\delta^2 & 0 & 0 & 0\\
0 & 0 & 0 & 1+b_4\delta^2\\
0 & 0 & b_3\delta^2 & 0
\end{pmatrix},
$$
\begin{align*}
H_{12}(\delta)=-i\alpha_0^{-2}\delta^2P_H(0)\sigma_1Q_H(0),
\quad H_{21}(\delta)=-i\alpha_0^{-2}\delta^2Q_H(0)\sigma_1P_H(0),
\end{align*}
where
\begin{align*}
& b_1=\alpha_0^{-2}\theta_1\|Q_c\|_{L^2(\R)}^2, \quad
b_2=-\alpha_0^{-2}\theta_1\|\pd_cQ_c\|_{L^2(\R)}^2,\\
& b_3=-4\alpha_0^{-4},\quad
b_4=\alpha_0^{-2}\|sQ_c\|_{L^2(\R)}^2\|Q_c\|_{L^2(\R)}^{-2},
\quad \sigma_1=\begin{pmatrix}  0 & 1\\ 1 & 0\end{pmatrix}.
\end{align*}
\par

First, we investigate the spectrum of $H_{11}(\delta)$.
Suppose $\lambda$ is an eigenvalue of the matrix $H_{11}(\delta)$.
Then
\begin{align*}
& \det(\lambda I-H_{11}(\delta))\\ =&
\left\{(\lambda+2i\alpha_0^{-2}\delta)^2-b_1\delta^2-b_1b_2\delta^4 \right\}
\left\{(\lambda+2i\alpha_0^{-2}\delta)^2-b_3\delta^2-b_3b_4\delta^4\right\}=0.
\end{align*}
Hence there exist eigenvalues $\hat{\lambda}_i$ $(i=1,2,3,4)$ of
$H_{11}(\delta)$ satisfying
\begin{align*}
& \hat{\lambda}_1
=-\delta\left(2i\alpha_0^{-2}-\alpha_0^{-1}\gamma+O(\delta^2)\right),
\quad
\hat{\lambda}_2=
-\delta\left(2i\alpha_0^{-2}+\alpha_0^{-1}\gamma+O(\delta^2)\right),
\\ &
\hat{\lambda}_3=-4i\alpha_0^{-2}\delta\left(1+O(\delta^2)\right),\quad
\hat{\lambda}_4=O(\delta^3).
\end{align*}
\par

Let $R_{ii}(\lambda,\delta)=(\lambda-H_{ii}(\delta))^{-1}$ for $i=1,2$
and let
\begin{gather*}
R_0(\lambda,\delta)=
\begin{pmatrix}
R_{11}(\lambda,\delta) & 0\\
0 & R_{22}(\lambda,\delta)
\end{pmatrix},\\
V_0(\lambda,\delta)=
\begin{pmatrix}
0 & H_{12}(\lambda,\delta)R_{22}(\lambda,\delta)\\
H_{21}(\lambda,\delta)R_{11}(\lambda,\delta) & 0
\end{pmatrix}. 
\end{gather*}
We remark that $R_{22}(\lambda,\delta)$ is uniformly bounded for
$\lambda\in U$ and $\delta\in(0,\delta_0)$.
Suppose that $|\lambda-\hat{\lambda_i}|=c_1\delta$, where $c_1\in
(0,\alpha_0^{-1}|\gamma|\delta/4)$ is a constant such that
$|\hat{\lambda}_j-\hat{\lambda}_k|\ge c_1\delta$
for every $j$, $k=1,2,3,4$ with $j\ne k$.
Then in view of the definitions of $H_{12}(\lambda,\delta)$ and
$H_{21}(\lambda,\delta)$,  we have
\begin{equation}
  \label{eq:5.6}
\left\|V_0(\lambda,\delta)\right\|_{B(L^2(\R))}=O(\delta),
\end{equation}
and
\begin{equation}
  \label{eq:5.7}
(\lambda-H(\delta))^{-1}=
R_0(\lambda,\delta)\sum_{i=0}^\infty V_0(\lambda,\delta)^i.
\end{equation}
Now let
\begin{align*}
P_{H,i}(\delta)=& \frac{1}{2\pi 
i}\oint_{|\lambda-\hat{\lambda}_i|=c_1\delta}
(\lambda-H(\delta))^{-1}d\lambda,
\\
\widehat{P}_{H,i}(\delta)=&
\frac{1}{2\pi i}\oint_{|\lambda-\hat{\lambda}_i|=c_1\delta}
R_0(\lambda,\delta)d\lambda.
\end{align*}
Combining \eqref{eq:5.6} and \eqref{eq:5.7} with
the fact that
\begin{align*}
\left\|R_0(\lambda,\delta)V_0(\lambda,\delta)\right\|_{B(L^2(\R))}
=&
\left\|
\begin{pmatrix}
0 & R_{11}H_{12}R_{22} \\
R_{22}H_{21}R_{11} & 0
\end{pmatrix}\right\|_{B(L^2(\R))}
=O(\delta),
\end{align*}
we have
$$
\left\| P_{H,i}(\delta)-\widehat{P}_{H,i}(\delta)\right\|=O(\delta)
\quad\text{for every $i= 1,2,3,4$.}
$$
Hence it follows that $R(\widehat{P}_{H,i}(\delta))$ is isomorphic to
$R(P_{H,i}(\delta))$ and that \linebreak
$R(P_{H,i}(\delta))$ is 1-dimensional for $i=1,2,3,4$.
Furthermore, we see that eigenvalues of $H(\delta)$ which lie in $U$
satisfy $|\lambda-\hat{\lambda}_i|< c_1\delta$ for an $i\in\N$ with
$1\le i\le 4$.

 Since $d\|Q_c\|_{L^2(\R)}^2/dc>0$ for $p\in(1,5)$, we see that
$\gamma$ is a positive number and that there exist eigenvalues $\lambda_1$ and $\lambda_2$ satisfying
$$
\alpha_0^{-1}\gamma\delta/2<\Re\lambda_1< 3\alpha_0^{-1}\gamma\delta/2,
\quad
-3\alpha_0^{-1}\gamma\delta/2<\Re\lambda_2< -\alpha_0^{-1}\gamma\delta/2.$$
Thus we complete the proof of Lemma \ref{lem:5.2}.
\end{proof}

\begin{proposition}
\label{prop:5.3}
Let $j$, $m\in\N$, $\varepsilon=m^{-1}$ and $\delta=j\varepsilon$.
Let $\beta=\min(p-1,1)/6$. Then there exists an
$m_*\in\N$ such that if $m\ge m_*$, the linearized operator
$\mathcal{H}(\varepsilon,\delta)$ with $j=[m^\beta]$ has an unstable
eigenvalue.
\end{proposition}

\begin{proof}
In order to prove Proposition \ref{prop:5.3}, we will show the
spectrum of $\mathcal{H}(\varepsilon,\delta)$ becomes close to the spectrum 
of
$H(\delta)$ as $\varepsilon\downarrow0$.  Let
$$
\mathcal{H}_0=i
\begin{pmatrix}
\frac{-2jm}{r^2} & \Delta_r-\omega-\frac{m^2+j^2}{r^2}\\
\Delta_r-\omega-\frac{m^2+j^2}{r^2} & \frac{-2jm}{r^2}
\end{pmatrix},
$$
and $H_0=U\mathcal{H}_0U^{-1}$.
Let
$$
\mathcal{D}(\lambda)=(\tau_{\bar{r}}\tilde{\chi}_0)(\lambda-H_0)^{-1}
(\tau_{\bar{r}}\chi_0)
+\tau_{\bar{r}}\tilde{\chi}_1(\lambda-H(\delta))^{-1}\chi_1
\tau_{-\bar{r}}.$$
Then we have
$$
\mathcal{D}(\lambda)U(\lambda-\mathcal{H}(\varepsilon,\delta))U^{-1}
=I+R_3+R_4,$$
where
\begin{align*}
& R_3=i(\tau_{\bar{r}})\tilde{\chi}_0(\lambda-H_0)^{-1}\left\{
\begin{pmatrix} 0 & [\pd_r^2,\tau_{\bar{r}}\chi_0]
 \\ [\pd_r^2,\tau_{\bar{r}}\chi_0] & 0\end{pmatrix}
-(\tau_{\bar{r}}\chi_0)\phi_\omega^{p-1} \begin{pmatrix}
0 & 1 \\ p & 0 \end{pmatrix}\right\}
\\ &
R_{4}=i\tau_{\bar{r}}\tilde{\chi}_1(\lambda-H(\delta))^{-1}\left\{
\begin{pmatrix} 0 & [\pd_r^2,\chi_1] \\ [\pd_r^2,\chi_1] & 0\end{pmatrix}
-\chi_1(R_{41}+R_{42})\right\}\tau_{-\bar{r}},
\\ &
R_{41}=
\begin{pmatrix}
\frac{-2\delta}{(\alpha_0+\varepsilon r)^2}+\frac{2\delta}{\alpha_0^2}
& -\frac{1+\delta^2-\frac14\varepsilon^2}{(\alpha_0+\varepsilon r)^2}
+\frac{1+\delta^2}{\alpha_0^2}\\
-\frac{1+\delta^2-\frac14\varepsilon^2}{(\alpha_0+\varepsilon r)^2}
+\frac{1+\delta^2}{\alpha_0^2} &
\frac{-2\delta}{(\alpha_0+\varepsilon r)^2}+\frac{2\delta}{\alpha_0^2}
\end{pmatrix},\\
& R_{42}=\begin{pmatrix}0 & f(\phi_\omega)-f(Q_c) \\
f'(\phi_\omega)-f'(Q_c) & 0\end{pmatrix}.
\end{align*}
We remark that
\begin{gather*}
|[\pd_r^2,\chi_i]\|_{B(L^2(\R),H^{-1}(\R))}=O(l^{-1})
\quad\text{for $i=0,1$,}\\
\|\chi_1R_{41}\|_{B(\Lr)}+\|R_{42}\|_{B(\Lr)}=O(\varepsilon^{6\beta}l).
\end{gather*}
We have
$$\sup_{\lambda\in\C,|\lambda|\le\omega/2}
\|(\lambda-\mathcal{H}_0)^{-1}\|_{B(H^{-2}_r(\R^2),L^2(\R^2))}<\infty,$$
since
$${}^tO\mathcal{H}_0O=i \begin{pmatrix}
\Delta_r-\omega-\frac{(m+j)^2}{r^2} & 0 \\
0 & -\Delta_r+\omega+\frac{(m-j)^2}{r^2} \end{pmatrix},$$
where
$$O=\frac{1}{\sqrt{2}}\begin{pmatrix} 1 & 1 \\ 1 & -1\end{pmatrix}.$$
Lemma \ref{lem:5.2} yields that for $\delta\in(0,\delta_0)$,
there exists a $c>0$ such that
$$
\left\|(\lambda-H(\delta))^{-1}\right\|_{B(\Lr)}\le C\delta^{-1}$$
for every $\lambda\in U$ with
$\min_{1\le i\le 4}|\lambda-\lambda_i(\delta)|\ge c\delta$ and that
$\Re(\lambda_1(\delta)-c\delta)>0$.
Let $l=\delta^{-3}$. Then it follows from the above that
\begin{align*}
& \|R_3\|_{B(\Lr)}=O(\delta^3+e^{-2\sqrt{c}\delta^{-3}}),\\
& \|R_4\|_{B(\Lr)}=O(\delta^2+\varepsilon^{6\beta}\delta^{-4}).  
\end{align*}
Put
\begin{align*}
\mathcal{P}_{\mathcal{H},1}(\varepsilon,\delta)=& \frac{1}{2\pi i}
\oint_{|\lambda-\lambda_1(\delta)|=c\delta}
(\lambda-\mathcal{H}(\varepsilon,\delta))^{-1}d\lambda,
\\
\mathcal{P}_{H,1}(\varepsilon,\delta)=&
U^{-1}\tau_{\bar{r}}\tilde{\chi}_1P_{H,1}(\delta)\chi_1\tau_{-\bar{r}}U.
\end{align*}
Making use of Cauchy's theorem and noting that
$\delta\sim\varepsilon^{\beta}$, we have
\begin{align*}
& \left\|\mathcal{P}_{\mathcal{H},1}(\varepsilon,\delta)
-\mathcal{P}_{H,1}(\varepsilon,\delta) \right\|_{B(\Lr)}
\\ =&
\frac{1}{2\pi}\left\|\oint_{|\lambda|=c\delta}
\left\{(\lambda-\mathcal{H}(\varepsilon,\delta))^{-1}
-U^{-1}\mathcal{D}(\lambda)U\right\}
d\lambda\right\|_{B(\Lr)}
\\ \le &
C\delta^{-1}\sup_{|\lambda|=c\delta}(\|R_3\|_{B(L^2(-\bar{r},\infty))}
+\|R_4\|_{B(L^2(-\bar{r},\infty))})
\\ \le & C(\delta+\varepsilon^{6\beta}\delta^{-5})
\\ =& O(\delta).
\end{align*}
From the above, we conclude that the range of
$\mathcal{P}_{\mathcal{H},1}(\varepsilon,\delta)$ is isomorphic to the
range  of $P_{H,1}(\delta)$ and that there exists an eigenvalue $\lambda$ of
$\mathcal{H}(\varepsilon,\delta)$ with $\Re\lambda>0$.
Thus we complete the proof of Proposition \ref{prop:5.3}.
\end{proof}
Now we are in position to prove Theorem \ref{thm2}.
\begin{proof}[Proof of Theorem \ref{thm2}]
Let $\mathfrak{L}$ be the linearized operator of
\eqref{eq:NLS} around $e^{i(\omega t+m\theta)}\phi_\omega$. Then
$$
\mathfrak{L}=i \begin{pmatrix}
\Delta-\omega+\beta_1(r) & e^{2im \theta}\beta_2(r) \\
-e^{-2im \theta}\beta_2(r) & -\Delta+\omega-\beta_1(r)
\end{pmatrix}.$$
Proposition \ref{prop:5.3} tells us that $\mathfrak L$ has unstable
eigenvalues if  $m\in\N$ is large and $p\in(1,5)$. On the other hand,
\cite{Mi3} tells us that $\mathfrak L$ has an unstable eigenvalue if
$p>3$. Hence it follows that $\mathfrak L$ has an unstable eigenvalue 
if $p>1$ and $m\in\N$ is sufficiently large.
\end{proof}

\section*{Acknowledgment}
The author would like to express his gratitude to Professor 
Jaeyoung Byeon for his useful advice.

\end{document}